\documentclass[12pt]{article}
\usepackage[cp866]{inputenc}
\usepackage[russian,english]{babel}
\usepackage{latexsym,amsfonts,amssymb,amsmath,amsfonts}
\usepackage{setspace}

 \newtheorem{thm}{Theorem}[section]
 
 \newtheorem{lem}[thm]{Lemma}

 \numberwithin{equation}{section}

\begin{document}

\begin{center}

{\bf Solvability of initial-boundary value problems for
non-autonomous evolution equations}

S.G. Pyatkov\footnote{ MSC Primary 35K90; Secondary 47D06;
34G10

Keywords: operator-differential equation, Cauchy problem,
non-autonomous evolution equation, maximal regularity,
initial-boundary value problem}

\bigskip

\parbox{10cm}{\small
 Abstract.  The initial-boundary value problems  for
linear non-autonomous first order evolution equations are
examined. Our assumptions provide a unified treatment which is
applicable to many situations, where the domains of the
operators may change with time. We study existence, uniqueness
and maximal regularity of solutions in Sobolev spaces. In
contrast to the previous results we use only the continuity
assumption on the operators in the main part of the equation.}

\end{center}

\section{Introduction}

Let $\{A(t)\}_{t\in [0,T]}$ be a family of closed linear
operators in a Banach space $X$. We consider the Cauchy problem
\begin{equation}\label{e1}
L(t)u=u_{t}-A(t)u -B(t)u=f,\
\end{equation}
\begin{equation}\label{e2}
u(0)=u_{0},
\end{equation}
where the family of operators $B(t):X\to X$, $t\in [0,T]$ is
subordinate in a certain sense to the family $A(t)$.  The most
known approach to the study of this problem is proposed in the
articles by Acquistapace P. and Terreni B. \cite{acq,acq1}.
Their approach goes back to the operator sum method of Da Prato
and Grisvard \cite{gri}. The main assumptions on the operator
family $A(t)$  in these articles are the so-called
Acquistapace-Terreni conditions (see \cite{acq,acq1}) connected
with the behavior of the resolvent and the H\"{o}lder
continuity of the family $\{A(t)\}$. Further developments of
this method are exposed in \cite{hie,hie1,gio}. Some results
are also presented in \cite[Sect. 6.8]{tan}. Similar results
under other weaker conditions also with the use of the
H\"{o}lder continuity of the family $\{A(t)\}$ (in a certain
sense) and the property that the space $(D(A(t)),X)_{\theta,p}$
obtained by the real interpolation method (see \cite{trib}) is
independent of $t$ for some $\theta\in (0,1)$ are presented in
\cite[Ch. 3]{yagi}. A series of articles is connected with the
minimal conditions on the family $\{A(t)\}$ when only
continuity of this family is required. In this case it is
assumed that the domains of $A(t)$ are independent of $t$. We
can refer, for instance, to \cite{ama1,ama2,ulm,pru,but,are}.
Some of these results  are exposed in \cite[Ch.6]{lun}. We also
refer to the book \cite[Ch.4]{amann}, where the reader can find
relevant  results as well as the bibliography. The Hilbert
space results devoted to the problem \eqref{e1}, \eqref{e2} are
often based on the Lax-Milgram theorem and the study of the
corresponding sesquilinear forms (see
\cite{laa,laa1,qua,fack}).

 Our approach is similar to that
described  in \cite[Ch.4, Sect. 3]{amann} in which the problem
is reduced to an abstract initial-boundary value problem. This
approach (see, for instance, \cite{gre}) is often used in the
study of abstract boundary control problems (see
\cite{eng,eng1} and the bibliography therein). We present
conditions on the operators $A$, $B$, and the boundary operator
below which ensure solvability of the corresponding
initial-boundary value problem and the problem
\eqref{e1}-\eqref{e2} as well under the minimal smoothness
assumptions on $A,B$ and the boundary operator.  We do not
require any H\"{o}lder continuity assumptions for the operator
$A$.

\section{Preliminaries}

Let $X, Y$ be  Banach spaces.  The symbol $L(X,Y)$  stands for
the space of linear continuous operators defined on  $X$ with
values in $Y$. If $X=Y$ then we use the notation  $L(X)$.  Let
$A:X\to X$ be a closed linear operator in $X$ with a dense
domain $D(A)$. The symbol $R(A)$ stands for the range of $A$.
Denote by $\sigma(A), \ \ \rho(A)$ the spectrum and the
resolvent set of $A$. Let ${\mathbb C}^{-}=\{z\in {\mathbb C}:\
Re\,z < 0\}$ ($C^{+}=\{z\in {\mathbb C}:\  Re\,z > 0\}$) and
let $\Sigma_{\theta}=\{z\in {\mathbb C}:\ |\arg\,z|<\theta\}$.

In what follows, we employ  the operators $A:X\to X$ ($X$ is a
Banach space) being the generators of analytic semigroups (see
\cite{lun}), in this case  we will  assume that $\rho(A)\supset
\overline{\Sigma_{\theta}}$ for some $\theta\geq \pi/2$ and
$$
\|\lambda(\lambda I - A)^{-1}\|_{L(X)} < M \ \forall  \lambda \in \overline{\Sigma_{\theta}},
$$
where $M>0$ is some constant  and $I$ is the identity.

Let $A:X\to X$ be a generator of an analytic semigroup. Put
$H_k=D(A^k)$ (the latter space is endowed with the graph norm).
We can also define the spaces $H_k$ for (see
\cite[Sect.5]{gri1}, \cite[Ch.6]{haase}) for $k<0$. The norm in
$H_{k}$ agrees with $\|(A-\lambda I)^{-k}u\|$, where $\lambda
\in \rho(A)$.  By the real interpolation method (see
\cite{trib,haase}) we can construct
$B^s_q=(H_m,H_k)_{\theta,q}$, with $1< q< \infty$, $k<s<m$, and
$\theta=\frac{m-s}{m+k}$ (see the properties of these spaces in
\cite[Sect. 5]{gri1}, Sect. 1.14, Sect. 1.15.4 in \cite{trib},
Prop. 1 in \cite{pya1}).

Define the space  $L_{p}(0,T;X)$ ($X$ is a Banach space) as the
space of strongly measurable functions, defined on  $[0,T]$
with values in  $X$ such that  $\int_0^T
\|u(t)\|_{X}^{p}dt<\infty$.  We use also the Sobolev spaces
$W_{p}^{s}(0,T;X)$ (see the definition, for instance, in
\cite{gri2,tri1}). The space of  continuous functions defined
on $[0,T]$ with values in $X$ is denoted by $C([0,T];X)$.

 A Banach space X is called a UMD space (the other names are $\zeta-convex$ and $HT$-spaces)
 if the Hilbert transform
 $Pf=\lim\limits_{\varepsilon\rightarrow
0}\int\limits_{|t-y|>\varepsilon}\frac{f(t)}{t-y}\,dt$  extends
to bounded operator on $L_{p}(R,X)$ for some (or equivalently,
for each) $p \in (1,\infty)$. All subspaces and quotient spaces
of $L_{q}(G, \mu)$ for $1 < q < \infty$ have the UMD property.
We can say that Sobolev spaces, Hardy spaces and other well
known spaces of analysis are UMD if they are reflexive.

 A collection  of operators $\tau\subset L(X,Y)$ ($X,Y$~are Banach
spaces)  is called $R$-bounded if there exists a constant $C_{p}$ such
that (see \cite{den2})
\[
\begin{array}{c}
\Bigl(\sum\limits_{\varepsilon_1,\varepsilon_2,\ldots,\varepsilon_N\in
\{-1,1\}} \Bigl\|\sum\limits_{j=1}^N \varepsilon_j T_j
x_j\Bigr\|^p\Bigr)^{\frac{1}{p}}\leq C_p
\Bigl(\sum\limits_{\varepsilon_1,\varepsilon_2,\ldots,\varepsilon_N\in
\{-1,1\}} \Bigl\|\sum\limits_{j=1}^N \varepsilon_j
x_j\Bigr\|^p\Bigr)^{\frac{1}{p}},\
\end{array}
\]
for all $N$, $T_{1},T_{2},\ldots,T_{N} \in \tau$  and
$x_{1},x_{2},\ldots, x_{N}\in X$. The least constant $C_{p}$ in
this inequality is denoted by $R(\tau)$ and is called the
$R$-bound of the family $\tau$ (see  equivalent
definitions in \cite{kun,den,den1,pru1}).
Note that this definition is independent of  $p$.

Next, we present some conventional results concerning with the
solvability of the Cauchy problem.

Consider the Cauchy problem
\begin{equation}\label{e4}
u_{t}-Lu=f, \ \  u|_{t=0}=u_{0}.
\end{equation}

We assume that $L$ is a generator of analytic semigroup and
that

(A) a family $\tau=\{\lambda(L-\lambda I)^{-1}:\ \lambda \in
\overline{\Sigma_{\theta_{0}}}\}$  is $R$-bounded for some
$\theta_{0}\geq \pi/2)$.

Denote the $R$-bound of this family by $M_{L}$.

\begin{thm} \label{th2}
 Let $X$ be a UMD space and let the condition {\rm (A)} hold.
Then, for every  $f\in L_{q}(0,T;X)$ and $u_{0}\in
B_{q}^{1-1/q}$, there exits a unique solution to the problem
\eqref{e4} such that $ u\in L_{q}(0,T;D(L))$, $u_{t}\in
L_{q}(0,T;X)$ and the estimate
$$
\|u_{t}\|_{L_{q}(0,T;X)}+\|Lu\|_{L_{q}(0,T;X)}\leq
C(\|f\|_{L_{q}(0,T;X)}+\|u_{0}\|_{B_{q}^{1-1/q}})
$$
holds. The constant $C$ depends on the constant $M_{L}$, $X$,
and $q$ and is bounded for bounded constants $M_{L}$.
\end{thm}
   The former claim results from  \cite[Theorem
3.2]{pru1}, \cite[Theorem 4.4]{den}. The last statement is
actually  follows from the known results  (see the proof of
Theorem 4.4 and the claim of Theorem 3.19 in \cite{den}).

In the following theorems we replace the problem \eqref{e4}
with the problem
\begin{equation}\label{e5}
u_{t}-Lu+\gamma u=f, \ \  u|_{t=0}=0, \ \ \gamma>0
\end{equation}
where $\gamma>0$ is a parameter and $L:X\to X$ is a generator
of an analytic semigroup.

\begin{thm}\label{th4}
 Let $X$ be a UMD space and let the condition {\rm (A)} hold.
Then, for every  $f\in L_{q}(0,T;X)$, there exits a unique
solution to the problem \eqref{e5} such that $ u\in
L_{q}(0,T;D(L))$, $u_{t}\in L_{q}(0,T;X)$ and the estimate
\begin{equation}\label{e7}
\|u_{t}\|_{L_{q}(0,T;X)}+\|Lu\|_{L_{q}(0,T;X)}+\gamma
\|u\|_{L_{q}(0,T;X)} \leq C\|f\|_{L_{q}(0,T;X)}
\end{equation}
holds, where the constant $C$ is independent of $\gamma$. It
depends of the constant $M_{L}$,  $q$, and the space $X$.
\end{thm}
    We consider the operator $L-\gamma I$ rather
than  $L$. In order to  prove the claim, we should estimate the
quantity $R\{\lambda(L-\lambda-\gamma)^{-1},\ \lambda \in
\Sigma_{\theta_{0}}\}$ and employ Theorem \ref{th2}). First, we
can say that $R\{i \xi(L-i\xi-\gamma)^{-1},\ \xi\in {\mathbb
R}\}\leq R\{(i\xi+\gamma)(L-i\xi-\gamma)^{-1},\ \xi\in {\mathbb
R}\}\leq 2M_{L}$ in view of Kahane's contraction principle (see
Remark 2.3 in \cite{den2} and Lemma 3.5 in \cite{den1}) and the
definition of $R$-boundedness. Next, we refer to the inequality
 $R\{\lambda(L-\lambda-\gamma)^{-1},\ Re\lambda\geq 0\}\leq R\{i
\xi(L-i\xi-\gamma)^{-1},\ \xi\in {\mathbb R}\}\leq 2M_{L}$
whose proof is presented  in Theorem 4.4 in \cite{den1}.

\section{Conditions on the data and auxiliary  results}

Now we can state the main conditions on the data of the problem
\eqref{e1}, \eqref{e2}.

First, we assume that there exists a Banach spaces $D\subset X$
and $Y$ and a family of linear operators $Q(t):D\to Y$  such
that

(1) $A(t)\in C([0,T];L(D,X)),\ \
Q(t)\in C([0,T];L(D,Y))$;

(2) the operators $A_{t}=A(t)|_{\ker Q(t)}):X\to X$ are the
generators of  analytic semigroups for every $t\in [0,T]$;

(3) $X$ is a UMD space and the family $\tau=\{\lambda
(-A_{t}+\lambda I)^{-1}: \ \lambda\in \overline{{\mathbb
C}^{+}}\}$ is $R$ bounded and $R(\tau)\leq M$, where the
constant $M$ is independent of $t\in [0,T]$.

Put $B_{q}^{s}=(D,X)_{1-s,q}$,
$H_{q}^{s,r}(\alpha,\beta)=W_{q}^{s}(\alpha,\beta;X)\cap
L_{q}(\alpha,\beta;B_{q}^{r})$. Introduce also the space
$H_{q}^{1,1}(0,T)$ of functions $u\in L_{q}(0,T;D)$ such that
$u_{t}\in L_{q}(0,T;X)$. Endow this space with the norm
$$
\|u\|_{H_{q}^{1,1}(0,T)}^{q}=\int_{0}^{T}\|u_{t}(t)\|_{X}^{q}+\|u(t)\|_{D}^{q}\,dt.
$$
The space  $H_{q}(a,T)$ ($a\in [0,T)$) consists  of functions
$u\in H_{q}^{1,1}(0,T)$ such that $u(t)=0$ for $t<a$ in the
case of $a>0$ and $u(0)=0$ for $a=0$. The norm in this space
coincide with that in $H_{q}^{1,1}(0,T)$. Define also the space
$H_{q}(a,b)$ ($0\leq a<b<T$) as the restriction of functions in
$u\in H_{q}(a,T)$ to the segment $[0,b]$. We endow the space
$H_{q}(a,b)$ with the norm $\inf\|\tilde{u}\|_{H_{q}(0,T)}$,
where the infimum is taken over all extensions $\tilde{u}$ of
$u\in H_{q}^{1,1}(a,b)$ to the whole segment $[0,T]$. Actually
this space consists of the functions  $u\in L_{q}(0,b;D)$ such
that $u_{t}\in L_{q}(0,b;X)$ and $u=0$ for $t<a$. One more
equivalent norm in this space is the norm
$$
\|u\|_{H_{q}(a,b)}=\Bigl(\int_{a}^{b}\|u_{t}(t)\|_{X}^{q}+\|u(t)\|_{D}^{q}\,dt\Bigr)^{1/q}.
$$
Given a function $u\in H_{0}(a,b)$, define its extension to the
segment $[0,T]$ as follows:
$$
P_{0}^{b}u(t)=\left\{\begin{array}{cl}
u(t),\ & t\in [0,b),\\
u(2b-t),\ & t\in [b,\min(2b-a,T)],\\
0,\ & t\in (2b-a,T]\ \textrm{if}\ 2b-a<T.
\end{array}
 \right.
$$
There is the obvious inequality
$$
\|u\|_{ H_{q}(a,b)}\leq \|P_{0}^{b}u\|_{H_{q}(0,T)}\leq
2^{1/q}\|u\|_{H_{q}(a,b)}
$$
which allows to say that the   norms  $\|u\|_{H_{q}(a,b)}$ and
$\inf \|\tilde{u}\|_{H_{q}(a,T)}$, where the infimum is taken
over the set $\{\tilde{u}\in H_{q}(a,T)\ :\ \tilde{u}=u\
\textrm{a.e. on}\ (a,b)\}$, are equivalent. Moreover, the
constants in the corresponding inequalities are independent of
$a,b$. In what follows we use the latter norm as the norm in  the space
$H_{q}(a,b)$.

The following conditions on the perturbations $B(t)$ are similar to those in \cite{ama1}:

(4) $B(t)\in L_{1}(0,T;L(D,X))$ and there exists a continuous
function $\beta(\xi):[0,+\infty)\to {\mathbb R}$ such that
$\beta(0)=0$ and $\|B(t)u(t)\|_{L_{q}(a,b;X)}\leq
\beta(b-a)\|u(t)\|_{H_{q}(a,b)}$ for all $u\in H_{q}(a,b)$ and
$0\leq a<b\leq T$.

Given a function $g(t)$,
define the function $g_{\varepsilon}(t)=\left\{\begin{array}{c} g(t-\varepsilon),\ t\in [\varepsilon, T]\\
0,\ t\in [0,\varepsilon). \end{array}\right.$, where
$\varepsilon\in (0,T)$.
Next, we expose some additional conditions on the mapping $Q$.
We assume that there exists a Banach space $Z\subset L_{q}(0,T;Y)$ such that

(5) The mappings $Q: u(t)\to Q(t)u(t)$, $Q_{\tau}:u(t)\to
Q(\tau)u(t)$ $(\tau\in [0,T])$ belong to the class
$L(H_{q}^{1,1}(0,T),Z)$, the norms $\|Q_{\tau}\|_{L(H_{q}^{1,1}(0,T),Z)}$ are uniformly bounded   and  the
mapping $Q_{\tau}$ is surjective for every $\tau\in [0,T]$;

(6) for every $\varepsilon>0$, there exists $\delta>0$ such
that
$$
\|(Q_{\tau_{1}}-Q_{\tau_{2}})u\|_{Z}\leq \varepsilon \|u\|_{H_{q}^{1,1}(0,T)}
$$
for all $u\in H_{q}(0,T)$ and  $\tau_{1},\tau_{2}\in [0,T]$
such that  $|\tau_{2}-\tau_{1}|<\delta$;
$$
\|(Q-Q_{\tau})u\|_{Z}\leq \varepsilon \|u\|_{H_{q}^{1,1}(0,T)}
$$
for all  $u\in H_{q}(\tau,T)$ and $\tau,b$ such that ${\rm
supp\,} u\subset [\tau,b]$, $0\leq \tau<b\leq T$,
$b-\tau<\delta$;
$$
((Q-Q_{0})v)_{\varepsilon_{0}}\in Z,\  \|((Q-Q_{0})u)_{\varepsilon_{0}}\|_{Z}\leq \varepsilon \|u(t)\|_{H_{q}^{1,1}(0,T)}
$$
for all  $v\in H_{q}^{1,1}(0,T)$, some $\varepsilon_{0}\in
(0,T)$, and every  $u\in H_{q}(0,T)$ such that ${\rm supp\,}
u\subset [0,b]$ with  $b<\delta$.

The conditions (6) is new in contrast to the conditions (1), (2) for the mapping $Q(t)$ which were used, for instance, in \cite{eng}.
They characterize the continuity of the mappings $Q, Q_{\tau}$.

Next, we specify  some additional function spaces and describe
their properties. Let $g(t)\in Z$.  Fix $\varepsilon\in (0,T)$
and define the space $Z_{q}(0,T)$ as the subspace of functions
$g\in Z$ such that there exists $\varepsilon>0$ such that
$g_{\varepsilon}\in Z$. Below we demonstrate that if
$g_{\varepsilon}\in Z$ for some $\varepsilon>0$ then
$g_{\varepsilon}\in Z$ for all $\varepsilon>0$. So it is
natural to fix $\varepsilon_{0}>0$ and  introduce the norm
$\|g(t)\|_{Z_{q}(0,T)}=\|g(t)\|_{Z}+\|g_{\varepsilon_{0}}(t)\|_{Z}$.
By $Z_{q}(a,T)$ $(a>0)$ we mean the subspace of $Z$  comprising
the functions vanishing for $t<a$ which is endowed with the norm of $Z$. At last we denote by $Z_{q}(a,b)$
($0\leq a<b\leq T$) the space of functions $g(t)$ such that
there exist a function $\tilde{g}\in Z_{q}(a,T)$ agreeing with
$g$ on $(a,b)$ almost everywhere (a.e.). We put
$$
\|g\|_{Z_{q}(a,b)}=\inf_{\tilde{g}=g\ on\  (a,b)}\|\tilde{g}\|_{Z_{q}(a,T)}.
$$
It is possible to define the operator $Q_{\tau}$ and $Q$ on the
space $H_{q}(a,b)$ putting
$Q_{\tau}u(t)=Q_{\tau}\tilde{u}(t)|_{(a,b)}$,
$Qu(t)=Q\tilde{u}(t)|_{(a,b)}$ where  $u\in H_{q}(a,b)$ and
$\tilde{u}$ is an extension of $u$ to the segment $[0,T]$.
Clearly, this definition is correct.

Let $L_{\tau}=\partial_{t}-A(\tau)$,
$L_{\tau,\gamma}=\partial_{t}-A(\tau)+\gamma I$ ($\gamma>0$).
Denote by $L_{\tau}^{-1}: L_{q}(a,b;X)\to L(H_{q}(a,b))$ the
operator taking a function $f\in L_{q}(a,b;X)$ into a solution
to the problem
\begin{equation}
\label {eq}
\partial_{t}u(t)-A_{\tau}u(t)=f,\ \  u(a)=0
\end{equation}
from the class $H_{q}(a,b)$.
 This operator is well-defined due to Theorem
\ref{th2} and the conditions (1)-(3). Similarly, we can define
the operator $L_{\tau,\gamma}^{-1}$. Introduce in $H_{q}(a,b)$
the equivalent norm
$$
\|u\|_{H_{q,\gamma}(a,b)}=\|u\|_{H_{q}(a,b)}+\gamma\|u\|_{L_{q}(a,b;X)}.
$$
The space $H_{q}(a,b)$ endowed with this norm is denoted by
$H_{q,\gamma}(a,b)$.

\begin{lem} \label{l0} Let the conditions {\rm (1)-(3)} hold.
Then   the  norms
$$
\|L_{\tau}\|_{L(H_{q}(a,b),L_{q}(a,b;X))},\ \|L_{\tau}^{-1}\|_{L(L_{q}(a,b;X),H_{q}(a,b))},\
\|L_{\tau,\gamma}^{-1}\|_{L(L_{q}(a,b;X),H_{q,\gamma}(a,b))}$$
 are uniformly bounded by some constant $C$ independent of $\tau\in
[0,T]$, $0\leq a<b\leq T$, and $\gamma$.
\end{lem}
   The first norms  are bounded  due to the condition
(1). Given $f\in L_{q}(a,b;X)$, let $\tilde{f}=f$ for $t\in
(a,b)$ and $\tilde{f}=0$ for $t\not\in (a,b)$. By Theorem
\ref{th2} we can find a function $u\in H_{q}(0,T)$ such that
$L_{\tau}u(t)=\tilde{f}$,\ \ $u(0)=0$, and $u(t)\in
D(A_{\tau})$ for almost all $t$. Due to uniqueness of solutions
to the Cauchy problem $u=0$ for $t\leq a$ and $u$ satisfies the
estimate
$$
\|u\|_{H_{q}(a,b)}\leq \|u\|_{H_{q}(0,T)}\leq
C\|\tilde{f}\|_{L_{q}(0,T;X)}=C\|f\|_{L_{q}(a,b;X)},
$$
where the constant $C$ is independent of $\tau$ (see Theorem
\ref{th2}).  So the norms
$\|L_{\tau}^{-1}\|_{L_{q}(a,b;X),L(H_{q}(a,b))}$ are uniformly
bounded. To prove the second statement, we can repeat the
arguments with the use of Theorem \ref{th4}.

\begin{lem} \label{l1} Let the conditions {\rm (1)-(3), (5), (6)} hold. Then

a) the spaces $Z_{q}(a,b)$ are Banach spaces for all values of
the parameters $a,b$;

b) for all $\tau, a,b\in [0,T]$  $(a<b)$,
$Q_{\tau}(H_{q}(a,b))=Z_{q}(a,b)$;

c) there exists a constant $C>0$ independent of $\tau,a,b\in
[0,T]$ $(a<b)$  such that
$\|Q_{\tau}\|_{L(H_{q}(a,b),Z_{q}(a,b))}\leq C$;

d) for every $\tau\in [0,T]$, the space $H_{q}(a,b)$ is the
direct sum of the subspaces $H_{\tau k}(a,b)=\{u\in H_{q}(a,b):
u\in\ker Q_{\tau}\}$ and $H_{\tau d}(a,b)=\{u\in H_{q}(a,b):\
L_{\tau}u=0\}$;

e) the mapping $Q_{\tau}$ is an isomorphism of $H_{\tau
d}(a,b)$ onto $Z_{q}(a,b)$;  for every $a_{0}>0$,   the norms
$\|Q_{\tau}^{-1}\|_{L(Z_{q}(a,b),H_{\tau d}(a,b))}$ are
uniformly bounded by a constant independent of $\tau,a,b\in
[0,T]$,
 such that  $a_{0}\leq a<b\leq T$; the norms
$\|Q_{0}^{-1}\|_{L(Z_{q}(0,b),H_{0 d}(0,b))}$ are uniformly
bounded by a constant independent of $b\in (0,T]$.

f) if $g(t),g_{\varepsilon_{0}}(t)\in Z$ then
$g_{\varepsilon_{1}}(t)\in Z$ for every $\varepsilon_{1} >0$,
and, for a given $a>0$, there exists a constant $C$ independent
of $\tau, c,b\in [0,T]$ $(a\leq c<b\leq T)$ such that
$\inf_{\tilde{g}=g\ on\ (c,b)}
\|\tilde{g}(t)\|_{Z}+\|\tilde{g}_{\varepsilon_{0}}(t)\|_{Z}\leq
C \inf_{\tilde{g}=g\ on\  (c,b)}\|\tilde{g}(t)\|_{Z}$ for all
$g\in Z_{q}(c,b)$;

g) the norms $\inf_{\tilde{g}=g\ on\ (0,b)}
\|\tilde{g}(t)\|_{Z}+\|\tilde{g}_{\varepsilon_{0}}(t)\|_{Z}$ in
$Z_{q}(0,b)$ for different parameters $\varepsilon_{0}$ are
equivalent.
\end{lem}

The proof of a) is almost obvious. Examine, for example, the
case of $a=0, b=T$. Let $g_{n}\in Z_{q}(0,T)$ be a Cauchy
sequence in $Z_{q}(0,T)$ endowed with the norm $\|g(t)\|_{Z_{q}(0,T)}=\|g(t)\|_{Z}+\|g_{\varepsilon_{0}}(t)\|_{Z}$. Since $Z$ is a Banach space, there
exist functions $g(t)\in Z$, $\tilde{g}_{\varepsilon_{0}}(t)\in
Z$ such that $g_{n}\to g$, $g_{n\varepsilon_{0}}\to
\tilde{g}_{\varepsilon_{0}}$ as $n\to \infty$ in $Z$. In the
space  $L_{q}(0,T;Y)$ this convergence also takes place.
Extracting a subsequence if necessary we can assume  that
$g_{n\varepsilon_{0}}(t)\to \tilde{g}_{\varepsilon_{0}}$ a.e.
in $Y$. Since $g_{n\varepsilon_{0}}=0$ a.e. for
$t<\varepsilon_{0}$, we have that $\tilde{g}_{\varepsilon_{0}}(t)=0$
a.e. for $t<\varepsilon_{0}$.
 Extracting one more subsequence if necessary we can assume that $g_{n}(t)\to g(t)$ a.e. in $Y$.
 In this case $g_{n\varepsilon_{0}}(t)=g_{n}(t-\varepsilon_{0})\to g(t-\varepsilon_{0})$ for a. a. $t\geq \varepsilon_{0}$ in $Y$.
Hence, $\tilde{g}_{\varepsilon_{0}}=g_{\varepsilon_{0}}$ a.e.
and thus $g\in Z_{q}(0,T)$. The remaining proofs in the case a)
employ the same arguments.

 Fix $\tau\in [0,T]$.
Let $u\in H_{q}(a,b)$, $a\geq 0$. Demonstrate that $Q_{\tau}u\in Z_{q}(a,b)$. There exists an extension
$\tilde{u}\in H_{q}(a,T)$ of this function. If $a>0$ then
$Q_{\tau}\tilde{u}=0$ a.e. on $(0,a)$ and $Q_{\tau}\tilde{u}\in
Z$. By definition, $Q_{\tau}\tilde{u}\in Z_{q}(a,T)$. If
$a=0$ then $\tilde{u}_{\varepsilon_{0}}\in H_{q}(\varepsilon_{0},T)$ and
thus $Q_{\tau}\tilde{u}_{\varepsilon_{0}}\in Z$,  i.~e.,
$Q_{\tau}\tilde{u}\in Z_{q}(0,T)$ and $Q_{\tau}u\in
Z_{q}(0,b)$. Proceed with c). Let $a>0$. In view of the condition (5), we infer
$
\|Q_{\tau}\tilde{u}\|_{Z}\leq C\|\tilde{u}\|_{H_{q}(a,T)},
$
where $\tilde{u}\in H_{q}(a,T)$ is an extension of $u$. Taking
the infimum over all $\tilde{u}$, we obtain the estimate
$
\|Q_{\tau}\|_{L(H_{q}(a,b),Z_{q}(a,b))}\leq C,
$
where the constant $C$ is independent of $\tau$ and $a,b$.
Let $a=0$. In this case we have the inequality
$$
\|Q_{\tau}\tilde{u}\|_{Z}+ \|Q_{\tau}\tilde{u}_{\varepsilon_{0}}\|_{Z} \leq C_{1}\|\tilde{u}\|_{H_{q}(0,T)},
$$
where the constant $C_{1}$ is independent of $\varepsilon_{0},b,$ and $\tau$. As above
we arrive at the estimate
$
\|Q_{\tau}\|_{L(H_{q}(0,b),Z_{q}(0,b))}\leq C_{1}.
$
We have proven c).

Demonstrate that the mapping $Q_{\tau}:H_{q}(a,b)\to
Z_{q}(a,b)$ is surjective.
 Let $g\in Z_{q}(a,b)$ with $a>0$.
In this case there exists an extension $\tilde{g}\in
Z_{q}(a,T)$ such that $\tilde{g}=g$ a.e. on $(a,b)$. In view of
(5), there exists a function $u\in H_{q}^{1,1}(0,T)$ such that
$Q_{\tau}u(t)=\tilde{g}(t)$. Since $g(t)=\tilde{g}(t)=0$ a.e.
on $(0,a)$, $u(t)\in ker\, Q_{\tau}$ a.e. for $t\leq a$  and
thus $u(t)\in D(A_{\tau})$ for a.a. $t\in (0,a)$. Therefore,
$u\in L_{q}(0,a;D(A_{\tau})), u_{t}\in L_{q}(0,a;X)$. By
Theorem III 4.10.2 in \cite{amann}, $u\in
C([0,a];(D(A_{\tau}),X)_{1/q,q})$  after a possible
modification on a zero measure set and there  exists the trace
 $u(a)\in (D(A_{\tau}),X)_{1/q,q}$. Consider the problem
\begin{equation}\label{e8}
v_{t}-A_{\tau}v=0, \ \  v|_{t=a}=u(a).
\end{equation}
In view of the condition (3),  Theorem \ref{th2} implies that
there exists a unique solution to this problem such that  $v\in
H_{q}^{1,1}(a,T)$ and $v\in D(A_{\tau})$ a.e. on $(a,T)$. Thus,
$v(t)\in ker \,Q_{\tau}$ for almost all $t$. By construction,
we have that $Q_{\tau}(u-v)=Q_{\tau}(u)=\tilde{g}(t)$. Define
the function $\tilde{u}(t)=(u-v)(t)$  for $t\geq a$ and
$\tilde{u}(t)=0$ for $t<a$. Obviously,  $\tilde{u}\in
H_{q}(a,T)$. Thus, we have found a function $\tilde{u}\in
H_{q}(a,T)$ such that $Q_{\tau}(\tilde{u})=g(t)$ and have
proven that $Q_{\tau}(H_{q}(a,b))=Z_{q}(a,b)$ for $a>0$.

Consider the case of $a=0$. Let $g\in Z_{q}(0,b)$. There exists
a function $\tilde{g}(t)\in Z$ such that $g(t)=\tilde{g}(t)$ on
$(0,b)$ a.e.
 In view of (5), there exist $\tilde{u}(t)\in H_{q}^{1,1}(0,T)$ such that
$Q_{\tau}\tilde{u}(t)=\tilde{g}$. As in the previous case, we
can find a function $\tilde{v}_{\varepsilon_{0}}\in
H_{q}(\varepsilon_{0},T)$ such that
$Q_{\tau}\tilde{v}_{\varepsilon_{0}}=\tilde{g}_{\varepsilon_{0}}$.
In particular, we have that
$\tilde{v}_{\varepsilon_{}}(\varepsilon_{0})=0$.  Find a
constant $\delta>0$ such that $\delta<T-\varepsilon_{0}$ and
construct a scalar function  $\varphi(t)\in C^{\infty}([0,T])$
such that $\varphi(t)=1$ for $t\in [0,\delta]$, $\varphi(t)=0$
for $t\in [T-\varepsilon_{0},T]$. Let
$u(t)=\varphi(t)\tilde{v}_{\varepsilon_{0}}(t+\varepsilon_{0})+(1-\varphi(t))\tilde{u}(t)\in
H_{q}(0,T)$. By construction,
$$
Q_{\tau}u(t)=\varphi(t)Q_{\tau}(\tilde{v}_{\varepsilon_{0}}(t+\varepsilon_{0}))+(1-\varphi(t))Q_{\tau}(\tilde{u}(t))=\tilde{g}(t).
$$
Thus, we have found a function $u\in H_{q}(0,T)$ such that
$Q_{\tau}u=\tilde{g}$ and, therefore,  the mapping $Q_{\tau}$
is surjective for $a=0$ as well. Moreover, since $u(0)=0$,
$u_{\varepsilon_{1}}\in H_{q}(\varepsilon_{1},T)$ and
$Q_{\tau}u_{\varepsilon_{1}}(t)=\tilde{g}_{\varepsilon_{1}}(t)\in
Z$ for every $\varepsilon_{1}>0$. The latter (together with the
closed graph theorem) means that the definition of the space
$Z_{q}(0,T)$ (and $Z_{q}(0,b)$ as well) is independent of
$\varepsilon_{0}$ and different norms depending on the
parameter $\varepsilon$ are equivalent, i.~e., we have proven
g).

 Consider the case d). First, we note that  the mapping $Q_{\tau}:H_{\tau d}(a,b)\to
Z_{q}(a,b)$ is one-to-one. Indeed, if $u\in ker\,Q_{\tau}$ then
$u(a)=0$,  $u(t)\in D(A_{\tau})$ a.e., and $L_{\tau}u=0$.   By
Theorem \ref{th2}, $u\equiv 0$. Let $u\in H_{q}(a,b)$. By
Theorem \ref{th2}, there exist a solution to the problem
$L_{\tau}v=L_{\tau}u$, $v(a)=0$ from the class $H_{\tau
k}(a,b)$. We have $v=L_{\tau}^{-1}L_{\tau}u$. The maps
$I-L_{\tau}^{-1}L_{\tau}, L_{\tau}^{-1}L_{\tau}$ are bounded
and they are projections onto the respective subspaces $H_{\tau
d}(a,b)$, $H_{\tau k}(a,b)$. Hence, $H_{q}(a,b)= H_{\tau
d}(a,b)+H_{\tau k}(a,b)$, where the sum is direct.

To prove e), it suffices to  consider the case of $b=T$. The
general case follows from this one in view of the definitions
of the norms in $H_{q}(a,b)$ and $Z_{q}(a,b)$. Fix $a>0$ and
show  that the family of norms
$\|Q_{\tau}^{-1}\|_{L(Z_{q}(c,T),H_{\tau d}(c,T))}$ is bounded
on the set $(\tau,c)$, $\tau\in [0,T],c\in[a,T)$. Take an
arbitrary $\tau_{0}\in [0,T]$. The operator
$Q_{\tau_{0}}:H_{\tau_{0} d}(a,T)\to Z_{q}(a,T)$ is bounded,
one-to-one and surjective. The open mapping theorem ensures
that there exists a bounded inverse (see Corollary 2.12 in
\cite{rud}). Let $\|Q_{\tau_{0}}^{-1}\|_{L(Z_{q}(a,T),H_{\tau
d}(a,T))}=c_{\tau_{0}}$. Consider the equation
$Q_{\tau_{1}}(Q_{\tau_{0}}^{-1}\tilde{g})=g$. It can be
rewritten as follows
\begin{equation}\label{e9}
(Q_{\tau_{1}}-Q_{\tau_{0}})Q_{\tau_{0}}^{-1}\tilde{g}+\tilde{g}=g.
\end{equation}
In view of the property (6), for every $\varepsilon>0 $, there exists $\delta>0$ such that
\begin{equation}\label{e10}
\|(Q_{\tau_{1}}-Q_{\tau_{0}})Q_{\tau_{0}}^{-1}\tilde{g}\|_{Z_{q}(a,T)}\leq \varepsilon
c_{\tau_{0}}\|\tilde{g}\|_{Z_{q}(a,T)}, \ \  |\tau_{1}-\tau|<\delta.
\end{equation}
We can choose $\varepsilon$ so that  $\varepsilon c_{\tau_{0}}\leq 1/2$.
In this case the  equation \eqref{e9}
is solvable.
So for a fixed $\tau_{0}$ we have found a neighborhood (in the relative topology on $[0,T]$)
of the form $\{\tau_{1}\in [0,T]:\ |\tau_{1}-\tau|< \delta\}$ in  which the equation \eqref{e9}
is solvable and
\begin{equation}\label{e11}
\|\tilde{g}\|_{Z_{q}(a,T)}\leq 2\|g\|_{Z_{q}(a,T)}.
\end{equation}
The definition yields
$Q_{\tau_{1}}^{-1}g=(I-L_{\tau_{1}}^{-1}L_{\tau_{1}})Q_{\tau_{0}}^{-1}\tilde{g}$.
By Lemma \ref{l0} the norms
$\|L_{\tau_{1}}\|_{L(H_{q}(a,T),L_{q}(a,T;X)}$ and
$\|L_{\tau_{1}}^{-1}\|_{L_{q}(a,T;X),L(H_{q}(a,T))}$ are
uniformly bounded by some constant $C_{1}$ independent of
$\tau_{1}\in [0,T]$. Thus, there exists a
constant $C$  such that $
\|Q_{\tau_{1}}^{-1}g\|_{H_{q}(a,T)}\leq C \|g\|_{Z_{q}(a,T)}$
for all $|\tau-\tau_{1}|<\delta$. Therefore,
 for every $\tau\in [0,T]$, there exists a
neighborhood $\{\tau_{1}\in [0,T]:\ |\tau_{1}-\tau|< \delta\}$
in which $\|Q_{\tau_{1}}^{-1}\|_{L(Z_{q}(a,T), H_{q}(a,T)}\leq
C$ for some constant $C>0$. We can cover the whole segment
$[0,T]$ by neighborhoods with these properties and find a
finite subcovering. Thus, there exists a constant $C>0$ such
that $\|Q_{\tau}^{-1}\|_{L(Z_{q}(a,T), H_{q}(a,T)}\leq C$ for
every $\tau\in [0,T]$. Next, consider the norms
$\|Q_{\tau}^{-1}\|_{L(Z_{q}(b,T), H_{q}(b,T)}$ with $b>a$. All
these spaces $Z_{q}(b,T)$ are endowed  with the same norm
$\|\cdot\|_{Z}$ and $Z_{q}(b,T)\subset Z_{q}(a,T)$. Hence,
$$\|Q_{\tau}^{-1}\|_{L(Z_{q}(b,T), H_{q}(b,T))}\leq
\|Q_{\tau}^{-1}\|_{L(Z_{q}(a,T), H_{q}(a,T))}\leq C.
$$
Therefore, the norms $\|Q_{\tau}^{-1}\|_{L(Z_{q}(b,T),
H_{q}(b,T))}$ are bounded by a constant independent of $\tau,
b$.

Next, we consider the case of $a=0$. We fix
$0<\varepsilon_{0}<T$ and endow the spaces $Z_{q}(0,T)$  with
the norm
$\|g\|_{Z}+\|g_{\varepsilon_{0}}\|_{Z}=\|g\|_{Z_{q}(0,T)}$. The
difference with the previous case is that we use this new norm
and the parameter $\tau=0$ is fixed. The remaining arguments
are the same. Repeating the above arguments, we can establish
that there exists a constant $C>0$ such that
$$
\|Q_{0}^{-1}\|_{L(Z_{q}(0,T), H_{q}(0,T))}\leq C.
$$

Next, we fix $a>0$  and demonstrate that there exists a
constant $c$ such that $\|g\|_{Z}+\|g_{\varepsilon}\|_{Z}\leq c
\|g\|_{Z}$ for all $g\in Z_{q}(b,T)$ and all $b\geq a$. Indeed,
there exists a function $u\in H_{q}(b,T)$ such that
$Q_{\tau}u=g$, $u\in ker\,L_{\tau}$, $u(b)=0$. In this case the
function $u_{\varepsilon}\in H_{q}(b+\varepsilon,T)$ if
$b+\varepsilon<T$ and $u_{\varepsilon}=0$ otherwise. In any
case $g_{\varepsilon}=Q_{\tau}u_{\varepsilon}\in Z$ and
$\|g_{\varepsilon}\|_{Z_{q}(b,T)}\leq
c\|u_{\varepsilon}\|_{H_{q}(b,T)}\leq c\|u\|_{H_{q}(b,T)}\leq
c_{1}\|g\|_{Z}$. In the last inequality we use the equality
$u=Q_{\tau}^{-1}g$ and the above estimate for the norms
$\|Q_{\tau}^{-1}\|_{L(Z_{q}(b,T),H_{q}(b,T))}$. In the case of
the spaces $Z_{q}(c,b)$, we can use the inequality obtained for
the extension $\tilde{g}$ of a given functions $g$ and taking
the infimum we infer
$$
\inf_{\tilde{g}=g\ on\  (c,b)} \|\tilde{g}\|_{Z}+\|\tilde{g}_{\varepsilon}\|_{Z}\leq c \inf_{\tilde{g}=g\ on\  (c,b)} \|\tilde{g}\|_{Z}.
$$

\begin{lem} \label{l2} Assume that the conditions {\rm (1)-(3), (5), (6)} hold.
Then $Q\in L(H_{q}(a,d),Z_{q}(a,d))$ and there exists a
constant $c$ independent of $a,d$ with $0\leq a<d\leq T$ such
that $\|Q\|_{L(H_{q}(a,d)), Z_{q}(a,d))}\leq c$.
\end{lem}

   Let $u\in H_{q}(a,d)$ and  $a>0$.
There exists an extension $\tilde{u}\in H_{q}(a,T)$.
 In view of (5)  $Q\tilde{u}\in Z$ and $Q\tilde{u}=0$ for $t<a$.
By definition, $Q\tilde{u}\in Z_{q}(a,T)$. Moreover, the
condition (5) implies that $\|Q\tilde{u}\|_{Z_{q}(a,T)}\leq
c\|\tilde{u}\|_{H_{q}(a,T)}$, where the constant   $c>0$ is
independent of $a>0$. Taking the infimum over all $\tilde{u}$
we obtain the claim.
Let $a=0$. By Lemma \ref{l1} b),
$Q_{0}\tilde{u}\in Z_{q}(0,T)$. We have that
$(Q-Q_{0})\tilde{u}\in Z$ and
$((Q-Q_{0})\tilde{u})_{\varepsilon_{0}}\in Z$ in view of (6).
Thereby,   $(Q-Q_{0})\tilde{u}\in Z_{q}(0,T)$ and thus
$Q\tilde{u}\in Z_{q}(0,T)$. Show the estimate from the claim of
the lemma. First, in view of (5)
\begin{equation}\label{e12}
\|Q\|_{L(H_{q}(0,T),Z)}+\|Q_{0}\|_{L(H_{q}(0,T),Z)}\leq c_{1}<\infty.
\end{equation}
Next, we employ (6).
Given $\varepsilon=1$, find $\delta>0$ such that
\begin{equation}\label{e13}
\|((Q-Q_{0})\tilde{u})_{\varepsilon_{0}}\|_{Z}\leq \|\tilde{u}(t)\|_{H_{q}(0,T)}
\end{equation}
for all $ \tilde{u}\in H_{q}(0,T): \
\textrm{supp\,}\tilde{u}\in [0,b]$, $ b<\delta\leq T$. Fix a
parameter $b<\delta$ and construct a scalar function
$\varphi(t)\in C^{\infty}[0,T]$ such that $\varphi(t)=1$ for
$t\in [0,b/2]$, $\varphi(t)=0$ for $t\in [3b/4,T]$. We have
\begin{equation}\label{e14}
Q\tilde{u}=(Q-Q_{0})\varphi\tilde{u}+Q_{0}\varphi\tilde{u}+Q(1-\varphi)\tilde{u}
\end{equation}
In view of \eqref{e13}, for the first summand on the right-hand
side we have the  estimate
\begin{equation}\label{e15}
\|((Q-Q_{0})\varphi\tilde{u})_{\varepsilon_{0}}\|_{Z}\leq c_{1}
\|\varphi(t)\tilde{u}(t)\|_{H_{q}(0,T)}\leq c_{2}\|\tilde{u}(t)\|_{H_{q}(0,T)}.
\end{equation}
The second summand belongs to $Z_{q}(b/2,T)$ and in view of
Lemma \ref{l1} f)
\begin{equation}\label{e16}
\|(Q_{0}(1-\varphi)\tilde{u})_{\varepsilon_{0}}\|_{Z}\leq c_{3}
 \|Q_{0}(1-\varphi)\tilde{u}\|_{Z}\leq c_{4}\|\tilde{u}\|_{H_{q}(0,T)}.
\end{equation}
Lemma \ref{l1} e) implies that
\begin{equation}\label{e161}
\|(Q_{0}\varphi\tilde{u})_{\varepsilon_{0}}\|_{Z}\leq c_{5}\|(\varphi\tilde{u})_{\varepsilon_{0}}\|_{H_{q}(0,T)}
\leq c_{6}\|\tilde{u}\|_{H_{q}(0,T)}.
\end{equation}
The estimates \eqref{e12}-\eqref{e161} and the definition of
the norm in $Z_{q}(0,T)$ yield
$$
\|Q\tilde{u}\|_{Z_{q}(0,T)}\leq c_{6}\|\tilde{u}\|_{H_{q}(0,T)}.
$$
Next, taking the infimum over all extensions $\tilde{u}$ (in the case of $d<T$) we obtain the claim.
In the case of $d=T$, we have that $\tilde{u}=u$ and the arguments are the same.

\begin{lem}\label{l3} Let the conditions {\rm (1)-(6)} hold.
For every $\varepsilon>0$, there exists  a constant $\delta>0$
such that
$$
\|(A(t)-A(a))u\|_{L_{q}(a,b;X)}\leq \varepsilon \|u\|_{H_{q}(a,b)},
$$
$$
\|B(t)u\|_{L_{q}(a,b;X)}\leq \varepsilon \|u\|_{H_{q}(a,b)}
$$
for all $u\in H_{q}(a,b)$ and $a,b$ such that $0 \leq a <b\leq
T$, $b-a<\delta$. Fix $\varepsilon_{0}>0$ and endow the space
$Z_{q}(0,T)$ with the norm norm
$\|g\|_{Z_{q}(0,T)}=\|g\|_{Z}+\|g_{\varepsilon_{0}}\|_{Z}$. The
norm in $Z_{q}(0,b)$ is defined with the use of this norm.
Then, for a given $\varepsilon>0$, there exists a parameter
$\delta>0$ such that
$$
\|(Q-Q_{0})v_{0}\|_{Z_{q}(0,b)}\leq \varepsilon \|v_{0}\|_{H_{q}(0,b)},
$$
for all $v_{0}\in H_{q}(0,b)$ and $b<\delta$. Fix $\tau_{0}>0$
and endow the space $Z_{q}(a,T)$ $(a\geq \tau_{0})$ with the
norm agreeing with that in $Z$. The norm in $Z_{q}(a,b)$ is
defined with the use of this norm. Then for a given
$\varepsilon>0$, there exists $\delta>0$ such that
$$
\|(Q-Q_{a})v_{0}\|_{Z_{q}(a,b)}\leq \varepsilon \|v_{0}\|_{H_{q}(a,b)}
$$
for all $v_{0}\in H_{q}(a,b)$ and $a,b$ such that $\tau_{0}\leq
a<b\leq T$, $b-a<\delta$.
\end{lem}
   The first statement is obvious due to the uniform
continuity of the family $A(t)$ and the conditions (1), (4). In
view of (6), for a given $\varepsilon>0$, there exists a
parameter $\delta>0$ such that
$$
\|(Q-Q_{0})v_{0}\|_{Z_{q}(0,b)}\leq \|(Q-Q_{0})P_{0}^{b}v_{0}\|_{Z_{q}(0,T)}\leq
\frac{\varepsilon}{2} \|P_{0}^{b}v_{0}\|_{H_{q}(0,T)}\leq \varepsilon \|v_{0}\|_{H_{q}(0,b)},
$$
for all $v_{0}\in H_{q}(0,b)$ and $b<\delta/2$. The last
statement is also obvious due to the condition (6).

\section{Main results}

First, we consider an  initial-boundary value problem.
\begin{equation}\label{e17}
u_{t}=A(t)u +B(t)u+f,
\end{equation}
\begin{equation}\label{e18}
Q(t)u(t)=g(t),\ \  u(0)=u_{0}.
\end{equation}
Clearly, the problem has no solutions for arbitrary data
$g,u_{0}$. So we have the natural consistency condition
\begin{equation}\label{e19}
g(t)-Qv(t)\in Z_{q}(0,T),\ \ g(t)\in Z,\ \
\end{equation}
where $v(t)\in H_{q}^{1,1}(0,T)$ is an arbitrary function such
that $v(0)=u_{0}$. We assume that
\begin{equation}\label{e20}
 u_{0}\in B_{q}^{1-1/q}=(D,X)_{1/q,q}.
\end{equation}
In this case there exists a function $v\in H_{q}^{1,1}(0,T)$
such that $v(0)=u_{0}$ (Theorem 1.8.3 in \cite{trib}). Note
that the condition \eqref{e19} does not depend on this function
$v$. Indeed, let $v_{1}\in  H_{q}^{1,1}(0,T)$ be one more
function such that $v_{1}(0)=u_{0}$. In this case
$g(t)-Qv_{1}(t)=g(t)-Qv - Q(v_{1}-v)$. The first summand
belongs to $Z_{q}(0,T)$ in view of \eqref{e19} and Lemma
\ref{l2} ensures that  $Q(v_{1}-v)\in Z_{q}(0,T)$. Thus,
$g(t)-Qv_{1}(t)\in Z_{q}(0,T)$. Moreover, note that the
condition \eqref{e19} is equivalent to the condition
\begin{equation}\label{e191}
g(t)-Q_{0}v(t)\in Z_{q}(0,T),\ \ g(t)\in Z,\ \
\end{equation}
Indeed, $g(t)-Qv(t)=g(t)-Q_{0}v(t)-(Q-Q_{0})v$. However,
$(Q-Q_{0})v\in Z_{q}(0,T)$ in view of (6) and the definitions.
So the conditions \eqref{e191} and \eqref{e19}  are equivalent.

\begin{thm}\label{th5}
Assume that $f\in L_{q}(0,T;X)$ and the conditions {\rm
\eqref{e19}, \eqref{e20}, (1)-(6)} hold.  Then  there exists a
unique solution $u\in H_{q}^{1,1}(0,T)$ to the problem
\eqref{e17}-\eqref{e19}. A solution satisfies the estimate
\begin{equation}\label{e21}
\|u\|_{H_{q}^{1,1}(0,T)}\leq c(\|g-Qv\|_{Z_{q}(0,T)}+\|u_{0}\|_{B_{q}^{1-1/q}}+ \|v\|_{H_{q}^{1,1}(0,T)}+\|f\|_{L_{q}(0,T;X)}),
\end{equation}
where the constant $c$ is independent of $g,u_{0}$, and  $f$,
and $v\in H_{q}^{1,1}(0,T)$ is an arbitrary function such that
$v(0)=u_{0}$.
\end{thm}

Construct a function  $v\in H_{q}^{1,1}(0,T)$ such that
$v(0)=u_{0}$ and make the change of variables $u=\omega+v$. The
function $\omega$ is a solution to the problem
\begin{equation}\label{e22}
\omega_{t}=A(t)\omega+B(t)\omega+f_{0},\ \ f_{0}=f-L(t)v,
\end{equation}
\begin{equation}\label{e23}
Q(t)\omega(t)=g(t)-Qv(t)=g_{0}(t),\ \  \omega(0)=0.
\end{equation}
We first prove  solvability of this problem on a small time
interval $[0,\gamma_{0}]$ and then extend a solution to the
whole segment $[0,T]$.
 We look for a solution $\omega\in H_{q}(0,\gamma)$ to the problem \eqref{e21},
\eqref{e22} in the form
$\omega=Q_{0}^{-1}\tilde{g}+L_{0}^{-1}\tilde{f}$ with
$\tilde{f}\in L_{p}(0,\gamma;X)$ and $\tilde{g}\in
Z_{q}(0,\gamma)$. Inserting this representation in \eqref{e22},
\eqref{e23}, we infer (recall that $Q_{0}^{-1}\tilde{g}\in
ker\,L_{0}$, $L_{0}^{-1}\tilde{f}\in {ker\,}Q_{0}$)
\begin{equation}\label{e24}
\tilde{f} =(A(t)-A(0))Q_{0}^{-1}\tilde{g}+B(t)(Q_{0}^{-1}\tilde{g}+L_{0}^{-1}\tilde{f}) +(A(t)-A(0))L_{0}^{-1}\tilde{f} +
f_{0},
\end{equation}
\begin{equation}\label{e25}
\tilde{g}=g_{0}(t)-(Q-Q_{0})L_{0}^{-1}\tilde{f}-(Q-Q_{0})Q_{0}^{-1}\tilde{g}.
\end{equation}
Rewrite \eqref{e24}, \eqref{e25} in the form
\begin{equation}\label{e26}
\tilde{f} =S(\tilde{f},\tilde{g})+f_{0},
\end{equation}
\begin{equation}\label{e27}
\tilde{g}=S_{0}(\tilde{f},\tilde{g}) +g_{0}.
\end{equation}
These equalities are a system of equation for recovering
 the unknowns $\tilde{f}\in
L_{q}(0,\gamma;X),\tilde{g}\in Z_{q}(0,\gamma)$. Study the
properties of operators on the right-hand side of these
equations. By Lemma \ref{l3}, for a given $\varepsilon>0$,
there exists a parameter $\delta>0$ such that
$$
\|(Q-Q_{0})v_{0}\|_{Z_{q}(0,\gamma)}\leq \varepsilon \|v_{0}\|_{H_{q}(0,\gamma)},
$$
for all $v_{0}\in H_{q}(0,\gamma)$ and $\gamma<\delta$. We take
$v_{0}=L_{0}^{-1}\tilde{f}+Q_{0}^{-1}\tilde{g}$. In this case
we arrive at the inequality (see Lemmas \ref{l0}, \ref{l1})
\begin{equation}\label{e28}
\|S_{0}(\tilde{f},\tilde{g})\|_{Z_{q}(0,\gamma)}\leq  C_{1}\varepsilon (\|\tilde{f}\|_{L_{q}(0,\gamma;X)}+\|\tilde{g}\|_{Z_{q}(0,\gamma)}),
\end{equation}
where the constant $C_{1}$ is independent of $\gamma$ and
$\gamma<\delta$. Next, we consider the operator
$S(\tilde{f},\tilde{g})$. Similarly, by Lemma \ref{l3}, for a
given  $\varepsilon_{1}>0$, there exists a parameter
$\delta_{1}>0$ such that
$$
\|(A(t)-A(0))v_{0}\|_{L_{q}(0,\gamma;X)}\leq \varepsilon_{1}\|v_{0}\|_{H_{q}(0,\gamma)},
$$
$$
\|B(t)v_{0}\|_{L_{q}(0,\gamma;X)}\leq \varepsilon_{1}\|v_{0}\|_{H_{q}(0,\gamma)},
$$
for $\gamma<\delta_{1}$. The definition of the operator $S$
yields
\begin{equation}\label{e29}
\|S(\tilde{f},\tilde{g})\|_{L_{q}(0,\gamma;X)}\leq  \varepsilon_{1}C_{2}(\|\tilde{f}\|_{L_{q}(0,\gamma;X)}+\|\tilde{g}\|_{Z_{q}(0,\gamma)}),
\end{equation}
In view of \eqref{e28}, \eqref{e29} the operator $R:
(\tilde{f},\tilde{g})\to (S(\tilde{f},\tilde{g}),
S_{0}(\tilde{f},\tilde{g})$ taking the space
$X_{\gamma}=L_{q}(0,\gamma;X)\times Z_{q}(0,\gamma)$ into
itself satisfies the estimate
\begin{equation}\label{e30}
\|R\|_{L(X_{\gamma},X_{\gamma})}\leq (\varepsilon C_{1}+\varepsilon_{1}C_{2})(\|\tilde{f}\|_{L_{q}(0,\gamma;X)}+\|\tilde{g}\|_{Z_{q}(0,\gamma)}).
\end{equation}
Thus, if we choose $\varepsilon,\varepsilon_{1}$ so that
$\varepsilon=1/4C_{1}$, $\varepsilon_{1}=1/4C_{2}$ then, for
$\gamma<\min(\delta,\delta_{1}) $,  the operator
$R:X_{\gamma}\to X_{\gamma}$ is a contraction  and thereby
there exists a unique solution $(\tilde{f},\tilde{g})$ to the
system \eqref{e26}, \eqref{e27} from the space $X_{\gamma}$. In
this case the corresponding function
$\omega=Q_{0}^{-1}\tilde{g}+L_{0}^{-1}\tilde{f}$ is a solution
to the problem \eqref{e22}, \eqref{e23} defined on the segment
$[0,\gamma]$. Fix this parameter $\gamma$ and denote it by
$\gamma_{0}$.

Next, we prove that there exists a number $\tau>0$ such that if
the problem \eqref{e22}, \eqref{e23} is solvable on the segment
on $[0,\gamma]$ ($\gamma\geq \gamma_{0}$) then it is solvable
on $[0,\min(T,\gamma+\tau)]$. Indeed, let $\omega\in
H_{q}(0,\gamma)$ be a solution to the problem \eqref{e22},
\eqref{e23}. Define a function
$\omega_{0}=P_{0}^{\gamma}\omega\in H_{q}(0,T)$ and make the
change of variables $\omega=\omega_{1}+\omega_{0}$. In this
case the function $\omega_{1}$ is a solution to the problem
\begin{equation}\label{e31}
\omega_{1t}=A(t)\omega_{1}+B(t)\omega_{1}+f_{1}\in L_{q}(0,T;X),\  f_{1}=f_{0}-L(t)\omega_{0},
\end{equation}
\begin{equation}\label{e32}
Q(t)\omega_{1}(t)=g_{0}(t)-Q(t)\omega_{0}=g_{1}(t)\in Z_{q}(\gamma,T),\ \  \omega_{1}(\gamma)=0.
\end{equation}
We now repeat the previous arguments.
 We look for a solution $\omega_{1}\in H_{q}(\gamma,\gamma_{1})$ to the problem \eqref{e31},
\eqref{e32} in the form
$\omega_{1}=Q_{\gamma}^{-1}\tilde{g}+L_{\gamma}^{-1}\tilde{f}$
with $\tilde{f}\in L_{p}(\gamma,\gamma_{1};X)$ and
$\tilde{g}\in Z_{q}(\gamma,\gamma_{1})$. Inserting this
representation in \eqref{e31}, \eqref{e32}, we infer (recall
that $Q_{\gamma}^{-1}\tilde{g}\in ker\,L_{\gamma}$,
$L_{\gamma}^{-1}\tilde{f}\in {ker\,}Q_{\gamma}$)
\begin{equation}\label{e33}
\tilde{f} =(A(t)-A(\gamma))Q_{\gamma}^{-1}\tilde{g}+B(t)(Q_{\gamma}^{-1}\tilde{g}+L_{\gamma}^{-1}\tilde{f}) +(A(t)-A(\gamma))L_{\gamma}^{-1}\tilde{f} +
f_{1},
\end{equation}
\begin{equation}\label{e34}
\tilde{g}=g_{1}(t)-(Q-Q_{\gamma})L_{\gamma}^{-1}\tilde{f}-(Q-Q_{\gamma})Q_{\gamma}^{-1}\tilde{g}.
\end{equation}
Rewrite \eqref{e33}, \eqref{e34} in the form
\begin{equation}\label{e35}
\tilde{f}=S(\tilde{f},\tilde{g})+f_{1},
\end{equation}
\begin{equation}\label{e38}
\tilde{g}=S_{0}(\tilde{f},\tilde{g}) +g_{1}.
\end{equation}
By Lemma \ref{l3}, for a given $\varepsilon>0$, there exists a
parameter $\delta>0$ such that
$$
\|(Q-Q_{\gamma})v_{0}\|_{Z_{q}(\gamma,\gamma_{1})}\leq \varepsilon \|v_{0}\|_{H_{q}(\gamma,\gamma_{1})},
$$
for all $v_{0}\in H_{q}(\gamma,\gamma_{1})$ and $\gamma_{0}\leq
\gamma<\gamma_{1}<\delta+\gamma$. We take
$v_{0}=L_{\gamma}^{-1}\tilde{f}+Q_{\gamma}^{-1}\tilde{g}$. In
this case we arrive at the inequality (see Lemmas \ref{l0},
\ref{l1})
\begin{equation}\label{e39}
\|S_{0}(\tilde{f},\tilde{g})\|_{Z_{q}(0,\gamma)}\leq  C_{1}\varepsilon (\|\tilde{f}\|_{L_{q}(0,\gamma;X)}+\|\tilde{g}\|_{Z_{q}(0,\gamma)}),
\end{equation}
where the constant $C_{1}$ is independent of $\gamma\geq
\gamma_{0}$ and $\gamma_{1}<\delta+\gamma$. Next, we consider
the operator $S(\tilde{f},\tilde{g})$. Similarly, by Lemma
\ref{l3}, for a given  $\varepsilon_{1}>0$, there exists a
parameter $\delta_{1}>0$ such that
$$
\|(A(t)-A(\gamma))v_{0}\|_{L_{q}(\gamma,\gamma_{1};X)}\leq \varepsilon_{1}\|v_{0}\|_{H_{q}(\gamma,\gamma_{1})},
$$
$$
\|B(t)v_{0}\|_{L_{q}(\gamma,\gamma_{1};X)}\leq \varepsilon_{1}\|v_{0}\|_{H_{q}(\gamma,\gamma_{1})},
$$
for $\gamma_{1}<\delta_{1}+\gamma$. The definition of the
operator $S$ yields
\begin{equation}\label{e40}
\|S(\tilde{f},\tilde{g})\|_{L_{q}(\gamma,\gamma_{1};X)}\leq  \varepsilon_{1}C_{2}(\|\tilde{f}\|_{L_{q}(\gamma,\gamma_{1};X)}+\|\tilde{g}\|_{Z_{q}(\gamma,\gamma_{1})}),
\end{equation}
In view of \eqref{e39}, \eqref{e40} the operator $R:
(\tilde{f},\tilde{g})\to (S(\tilde{f},\tilde{g}),
S_{0}(\tilde{f},\tilde{g})$ taking the space
$X_{\gamma,\gamma_{1}}=L_{q}(\gamma,\gamma_{1};X)\times
Z_{q}(\gamma,\gamma_{1})$ into itself satisfies the estimate
\begin{equation}\label{e41}
\|R\|_{L(X_{\gamma,\gamma_{1}},X_{\gamma,\gamma_{1}})}\leq
(\varepsilon C_{1}+\varepsilon_{1}C_{2})(\|\tilde{f}\|_{L_{q}(\gamma,\gamma_{1};X)}+\|\tilde{g}\|_{Z_{q}(\gamma,\gamma_{1})}).
\end{equation}
Thus, if we choose $\varepsilon,\varepsilon_{1}$ so that
$\varepsilon=1/4C_{1}$, $\varepsilon_{1}=1/4C_{2}$ then, for
$\gamma_{1}<\min(\delta,\delta_{1})+\gamma$ ($\gamma_{1}\leq
T$), the operator $R:X_{\gamma,\gamma_{1}}\to
X_{\gamma,\gamma_{1}}$ is a a contraction  and thereby there
exists a unique solution $(\tilde{f},\tilde{g})$ to the system
\eqref{e33}, \eqref{e34} from the space
$X_{\gamma,\gamma_{1}}$. In this case the corresponding
function
$\omega_{1}=Q_{\gamma}^{-1}\tilde{g}+L_{\gamma}^{-1}\tilde{f}$
is a solution to the problem \eqref{e31}, \eqref{e32} defined
of the segment $[\gamma,\gamma_{1}]$. Fix
$\gamma_{1}<\min(\delta,\delta_{1})+\gamma$ and put
$\tau=\gamma_{1}-\gamma$. We can see that the constant $\tau$
depends only on some absolute constants and is independent of
$\gamma\geq \gamma_{0}$.  Let $\omega=\omega_{0}+\omega_{1}$,
where we extend the function $\omega_{1}$ by zero for $t\leq
\gamma$. Now the existence of a solution to the problem
\eqref{e22}, \eqref{e23} on the whole segment $[0,T]$ results
from the above-proven.

Demonstrate the a solution is unique. Let $w$ be a solution to
the problem \eqref{e22}, \eqref{e23} with $f_{0}=0, g_{0}=0$.
By Lemma \ref{l1} d), this solution is representable as
$w=L_{0}^{-1}\tilde{f}+Q_{0}^{-1}\tilde{g}$ with $\tilde{f}\in
H_{1}(0,T)$ and $\tilde{g}\in Z_{q}(0,T)$. Inserting it into
\eqref{e22}, \eqref{e23}, we arrive at the system \eqref{e24},
\eqref{e25}. Repeating the above arguments we obtain that
$\tilde{f}=0$, $\tilde{g}=0$ on some segment $[0,\gamma]$ (the
segment on which a solution to the problem exists). Next,
repeating the arguments on the segments of the form
$[\gamma,\gamma+k\tau]$ with $\tau$ the above parameter and $k$
a positive integer we prove that $\omega\equiv 0$.

\begin{thm}\label{th6}
Assume that $f\in L_{q}(0,T;X)$,  $u_{0}\in
B_{0q}^{1-1/q}=(D(A_{0}),X)_{1/q,q}$, and the conditions {\rm
(1)-(6)} hold. Then  there exists a unique solution $u\in
H_{q}^{1,1}(0,T)$ to the problem \eqref{e1}-\eqref{e2} such
that $u(t)\in D(A(t))$ for a.a. $t\in [0,T]$. A solution
satisfies the estimate
\begin{equation}\label{e42}
\|u(t)\|_{H_{q}^{1,1}(0,T)}\leq c(\|u_{0}\|_{B_{0q}^{1-1/q}}+ \|f\|_{L_{q}(0,T;X)}),
\end{equation}
where the constant $c$ is independent of $g,u_{0}$.
\end{thm}
   We consider the problem \eqref{e17}, \eqref{e18}, where $g(t)=0$.
By Theorem \ref{th2}, there exists a solution to the problem
$L_{0}v=0$, $v(0)=u_{0}$ from the class $H_{q}^{1,1}(0,T)$ such
that $v\in ker\,Q_{0}$. This solution meets the estimate
\begin{equation}\label{e43}
\|v\|_{H_{q}^{1,1}(0,T)}\leq c\|u_{0}\|_{B_{0q}^{1-1/q}}.
\end{equation}
Note that there is the natural embedding $B_{0q}^{1-1/q}\subset
B_{q}^{1-1/q}$ and thus there exists a constant $c>0$ such that
\begin{equation}\label{e44}
\|u_{0}\|_{B_{q}^{1-1/q}}\leq c\|u_{0}\|_{B_{0q}^{1-1/q}}.
\end{equation}
For this function $v$ the condition \eqref{e191}  (and
respectively \eqref{e19}) of Theorem \ref{th5} is fulfilled
automatically. Applying Theorem \ref{th5}, we can find a
solution to the problem \eqref{e17}, \eqref{e18}. Since
$g(t)=0$, we conclude that $u(t)\in D(A(t))$ for a.a. $t\in
[0,T]$. Proceed with the estimate \eqref{e42}. The first
summand in the estimate \eqref{e21} is estimated as
\begin{equation}\label{e45}
\|Qv\|_{Z_{q}(0,T)}=\|(Q-Q_{0})v\|_{Z_{q}(0,T)}\leq c \|v\|_{H_{q}(0,T)}\leq \|u_{0}\|_{B_{0q}^{1-1/q}}.
\end{equation}
The remaining summands are estimated with the use of
\eqref{e43} and \eqref{e44} and thus the estimate \eqref{e42}
is proven.

{\bf Remark.} It is possible that condition (3) is not
fulfilled but the claims of Theorems \ref{th5}-\ref{th6} remain
valid (possibly in some  other spaces). We need the maximal
regularity property of the family $\{A(t)\}_{t\in [0,T]}$
(i.~e., the claim of Theorem \ref{th2}) holds for every of the
operators $A(t)$). If the family $\{A(t)\}_{t\in [0,T]}$
consists of generators of analytic semigroups then (see
\cite[Theorem 3.14 ]{gri}, \cite[Theorem 1]{pya1},
\cite[Theorem 2.7]{gri2},  in \cite[ Theorem 4]{uva}) this
family enjoys this property in the spaces
$\tilde{X}_{s}=B_{t,q}^{s}=(D(A(t)), X)_{1-s,q}$ and we can
 assume that these spaces are independent of $t$ for some $s$. So we can replace the condition (3) with
the condition

(3$'$) $\Sigma_{\theta_{0}}\subset \rho(A)$ for some
$\theta_{0}\geq \pi/2$, there exists a constant $c>0$ such that
$$
\|\lambda(\lambda I-A(t))^{-1}\|_{X}\leq c, \ \forall t\in [0,T],\ \lambda \in \Sigma_{\theta_{0}},
$$
and there exists  $s_{0}\in [0,1)$ such that the spaces $ B_{t
q}^{s_{0}}=(D(A(t)),X)_{1-s_{0}, q} $ coincide for all $t\in
[0,T]$.

Under the conditions (3$'$), we can reformulate  the conditions
(1), (2), (4)-(6) with the space $\tilde{X}_{s_{0}}$ rather
than $X$.

\section{Some applications}

The results of this section are very close to those in
\cite{den1}. They are not new.

We consider  vector-valued parabolic initial-boundary  value
problems of the form
\begin{gather}\label{th8}
u_{t}- A(t,x,D)u = f(t,x),\ \ x\in G\subset {\mathbb R}^{n},\  t\in (0,T),\\
 B_{j}(t, x,D)u = g_{j}(t,x)\ (j = 1, . . . ,m),\  x\in \Gamma=\partial G, \  t\in (0,T)\\
 u(0, x) = u_{0}(x), x \in G.
\end{gather}
 Here
 $G$  is a bounded domain in ${\mathbb R}^{n}$ with boundary $\Gamma\in C^{2m}$,
$A(t, x,D) = \sum_{|\alpha|\leq 2m}
a_{\alpha}(x,t)D^{\alpha}u$, $B_{j}(t, x,D) =\sum_{|\alpha|\leq
m_{j}}b_{j\alpha}(x,t) D^{\alpha}u(x,t)$, where $a_{\alpha}$
and $b_{j\alpha}$ are $L(E)$-valued variable coefficients and
$m_{j}<2m$.  Denote $S=(0,T)\times \Gamma$ and $Q=(0,T)\times
G$. Let $E$ be a UMD space. Put
$W_{q}^{s,r}(Q;E)=L_{q}(0,T;W_{q}^{r}(G,;E))\cap
W_{q}^{s}(0,T;L_{q}(G;E))$,
$W_{q}^{s,r}(S;E)=L_{q}(0,T;W_{q}^{r}(\Gamma;E))\cap
W_{q}^{s}(0,T;L_{q}(\Gamma;E))$.   Here and in what follows, we
use the conventional multi-index notation
$D^{\alpha}=\partial_{x_{1}}^{\alpha_{1}}\partial_{x_{2}}^{\alpha_{2}}\ldots
\partial_{x_{n}}^{\alpha_{n}}$ with
$\partial_{x_{1}}=\frac{\partial}{\partial x_{1}}$. We assume
the following conditions on the data:

(i) $f\in L_{q}(0,T;L_{q}(G;E))$ ($q\in (1,\infty)$),

(ii) $g_{j}\in W_{q}^{k_{j},2mk_{j}}(S;E)$,
$k_{j}=1-m_{j}/2m-1/2qm$,

(iii) $u_{0}\in W_{q}^{2m-2m/q}(G;E)$,

(iv) If $k_{j} > 1/q$ then $B_{j}(0, x)u_{0}(x) = g_{j}(0, x)$
for $x\in \Gamma$. If $k_{j} = 1/q$ then $
\int_{0}^{T}\|g_{j}(t,x)-B_{j}(0,x)v(t,x)\|_{L_{q}(\Gamma;E)}^{q}\,\frac{dt}{t}<\infty$,
where $v\in W_{q}^{1,2m}(Q;E)$ and $v(0,x)=u_{0}(x)$.

We start with the ellipticity assumptions. To this end, we
denote the principal part of $A$ by $A_{0}$,
$A_{0}(x,t,D)=\sum_{|\alpha|=2m} a_{\alpha}(x,t)D^{\alpha}$.
The unit outer normal to $\Gamma$ at $x\in \Gamma$
 is denoted by $\nu(x)$.
We use the following conditions:

 (v)  for all $t\in [0,T]$, $x\in \overline{G}$, and $\xi\in {\mathbb R}^{n}$, $|\xi| =
 1$, and some $\theta_{0}>\pi/2$
 we have $\sigma(A_{0}(x,t,\xi))\subset {\mathbb C}\setminus \Sigma_{\theta_{0}}$;

(vi) for all $t\in [0,T]$, $x\in \Gamma$, all $\xi\in {\mathbb
R}^{n}$ with $\nu(x)\cdot\xi=0$, all $\lambda \in
\Sigma_{\theta_{0}}$ and $h\in E^{m}$,  the ordinary
differential equation system in ${\mathbb R}^{+}=(0,+\infty)$
\begin{gather*}
 \lambda v(y) - A_{0}(t, x, \xi+\nu(x)\partial_{y})v=0, \ y
> 0,\  B_{j}(t, x, \xi+ \nu(x)\partial_{y})v(0) = h_{j},\
\end{gather*}
where $j = 1,2, \ldots, ,m,$ admits a unique solution $v\in
C([0,\infty);E)$ decreasing at infinity.

 Now we turn to smoothness assumptions on the
coefficients of $A$ and $B_{j}$. We assume that

(vii) there are $r_{k}, s_{k}\geq q$ with $1/s_{k} + n/2mr_{k}
< 1-k/2m$ such that $a_{\alpha}\in
L_{s_{k}}(0,T;L_{r_{k}}(G;L(E))$ for $|\alpha|=k < 2m$, and
$a_{\alpha}\in C(\overline{Q};L(E))$ for $|\alpha|= 2m$,
$b_{j\beta}\in C^{1-m_{j}/2m,2m-m_{j}}(\overline{S};L(E))$ for
$|\beta|\leq m_{j}$.

The following theorem holds.

\begin{thm} Under the conditions {\rm (i)-(vii)}, there exists a unique solution to the
problem \ref{th8} such that
$$
u\in L_{q}(0,T;W_{q}^{2m}(G;E)), \ \ u_{t}\in L_{p}(Q;E).
$$
\end{thm}

   We can refer to Theorem \ref{th5}. We need only
to  check the conditions of this theorem. Put $X=L_{q}(G;E)$,
$D=W_{q}^{2m}(G;E)$, $A(t)=A_{0}(x,t,D)$,
$B(t)=\sum_{|\alpha|<2m} a_{\alpha}(x,t)D^{\alpha}$,
$Y=\prod_{j=1}^{m}W_{q}^{2m k_{j}}(\Gamma;E)$,
$Z=\prod_{j=1}^{m}W_{q}^{k_{j},2m k_{j}}(S;E)$. The condition
(1) is obvious. Theorem 8.2 in \cite{den} ensures the
conditions (2), (3)   for the operator $A(t)-\lambda_{0}I$ with
sufficiently large $\lambda_{0}>0$. So without loss of
generality, we can assume that the condition (2) is fulfilled,
otherwise we make the change $u=e^{\lambda_{0}t}v$ to reduce
the problem to this case. The boundedness  of the constant
$R(\lambda(\lambda -A(t)),\ \lambda \in \overline{{\mathbb
C}^{+}})$ in dependence on $t\in [0,T]$ results from the
continuity of the property $R$-boundedness (see Proposition 4.2
in \cite{den}). The idea of the proof is presented in Sect. 7.3
in \cite{den}. To justify the condition (4), we can use Lemma
3.10 in \cite{den1}. The function $\beta$ in (4) is just a
function of the form $\beta(\xi)=c\xi^{\delta}$ with $c$ some
constant and $\delta$ is small parameter.  Actually, the
condition (4) is justified in \cite{den1}. However, the proof
there is not detailed. To prove, the condition (4), we should
use the Lemma 3.5 in \cite{den} and more or less conventional
arguments. To prove (5), we can also use Lemma 3.5 in
\cite{den}. All conditions in (6) result from the continuity
properties of the mapping $Q$ and the definition of the norm in
$Z$. The most difficult condition is the condition \eqref{e19}.
Demonstrate that the condition \eqref{e19} is fulfilled. Let
$v\in W_{q}^{1,2m}(Q)$ be such that $v(0,x)=u_{0}(x)$. Let
$\tilde{g}=(\tilde{g}_{1},\tilde{g}_{1},\ldots,
\tilde{g}_{m})$, $\tilde{g}_{i}=g_{i}-B_{i}(t,x)v|_{S}$. By
Prop. 5.11 in \cite{gri2} for $n=1$, $\tilde{g}(t)\in
Z_{q}(0,T)$ if and only if $\tilde{g}_{j}(0, x)=0$ if $k_{j} >
1/q$ and
$\int_{0}^{\delta_{0}}\|\tilde{g}_{j}(t,x)\|_{L_{q}(\Gamma;E)}^{q}\,\frac{dt}{t}<\infty$
for some $\delta_{0}>0$ if $k_{j} = 1/q$.

Acknowledgement. The author was supported by a grant for
scientific schools with young scientists participation of the
Yugra State University and by the Act 211 of the Government of
the Russian Federation, contract No. 02.A03.21.0011.

\bigskip
Yugra State University,  Chekhov st. 16, 628012,
Hanty-Mansiisk, Russia, Sobolev Institute of Mathematics,
Novosibirsk, Russia

email: pyatkov@math.nsc.ru, s\_pyatkov@ugrasu.ru


\begin{thebibliography}{23}

\bibitem{acq} P. Acquistapace  and  B. Terreni,  {\it A unified
    approach to abstract linear nonautonomous parabolic
    equations,} Rend. Sere. Mat. Univ. Padova, {\bf 78} (1987), 47-107.

\bibitem{acq1} P. Acquistapace, {\it Maximal regularity for
    abstratc   linear  non-autonomous parabolic equations,} J. of Funct. Anal. {\bf 60} (1985),
    168-210.

\bibitem{ama1} H. Amann, {\it Maximal regularity for
    nonautonomous     evolution equations}, Adv. Nonlinear Stud. {\bf 4} (2004), 417-430.

\bibitem{ama2} H. Amann {\it  Nonautonomous parabolic
    equations     involving measures},  J. of Math. Sciences. {\bf 30} (2005), 4780-4802.

\bibitem{amann} H. Amann, {\it Linear and quasilinear
    parabolic     problems}, Birkh\"{a}user Verlag,
    Basel-Boston-Berlin, {\bf 1}, 1995.

\bibitem{laa1} W. Arendt, D. Dier, H. Laasri, E. M. Ouhabaz,
    {\it Maximal Regularity for Evolution Equations Governed by
    Non-Autonomous Forms}, 2014.
    https://hal.archives-ouvertes.fr/hal-00797181v1.

\bibitem{are} W. Arendt, R. Chill, S. Fornaro,  C. Poupaud,
 {\it $L_{p}$-maximal regularity for non-autonomous evolution
 equations}, J. Diff. Equat., {\bf 237} (2007), 1-26.

\bibitem{but} A. Butti {\it On the Evolution Operator for a
Class of Non-autonomous Abstract Parabolic Equations}, J. of
Math. Anal. Appl., {\bf  170}  (1992), 115-137.


\bibitem{gri} G. Da Prato, P. Grisvard, {\it Sommes
    d'op\'{e}rateurs lin\'{e}iares et \'{e}quations
    diff\'{e}rentielles op\'{e}rationnelles}, J. Math.
    Pures Appl. {\bf 54} (1975) 305-387.

\bibitem{den} R. Denk, M. Hieber, J.
    Pr\"{u}ss, {\it R-boundedness, Fourier multipliers and problems
    of elliptic and parabolic type}, Mem. Amer. Math. Soc. {\bf
    166} (2003).

\bibitem{den1} R. Denk, M. Hieber, and J.  Pr\"{u}ss,
    {\it Optimal     $L_{p}-L_{q}$-estimates for parabolic boundary value problems with
    inhomogeneous data},  Math. Z.  {\bf 257} (2007), 93--224.

\bibitem{den2}
 {\it Denk R.,  Krainer T.} {\it $R$-boundedness, pseudodifferential
operators, and maximal regularity for some classes of partial
differential operators},~ Manuscripta Math., \textbf{124}
(2007), 319--342.


\bibitem{gio} D. Di Giorgio, A. Lunardi, and R. Schnaubelt
    {\it     Optimal regularity and Fredholm properties of abstract
    parabolic operators in $L_{p}$ spaces on the real line},
    Proc. of the London Math. Soc., {\bf 91} (2005), 703-737.


    \bibitem{eng}  K.-J. Engel, B. Kl\"{o}ss,  R. Nagel,  B.
    Fijav\v{z}, · E. Sikolya,     {\it  Maximal controllability for
    boundary control problems}, Appl. Math. Optim.,  {\bf 62}
    (2010), 205-227.

\bibitem{eng1}  K.-J. Engel,   B.   Fijav\v{z},     {\it
    Exact and positive controllability of boundary control
    systems},      Networks \& Heterogeneous Media, {\bf  12}(2) (2017), 319-337.

\bibitem{fack} S. Fackler, {\it J.-L. Lions' problem concerning
    maximal regularity of equations governed by non-autonomous
    forms}, Annales de l'Institut Henri Poincare (C) Non Linear
    Analysis, {\bf 34} (2017),  699-709.


\bibitem{gal} C. Gallarati and M. Veraar, {\it Maximal
regularity for non-autonomous equations with measurable
dependence on time}, Potential Analysis, {\bf  46} (2017),
 527-567.


\bibitem{gre} G. Greiner, {\it Perturbing the boundary
    conditions     of  a     generator},
 Houston J.  of Math. {\bf
    13} (1987), 213--229.


\bibitem{gri1} P. Grisvard, {\it Commutative de deux functeurs
    d'interpolation et applications},  J. Math.
    pures et appliq.,  {\bf 45} (1966), 143--206.

\bibitem{gri2} P. Grisvard, {\it Equations
    diffe\-ren\-ti\-el\-les abstraites},  Ann.
    Scient. Ec. Norm. Sup.   4$^{e}$ series, {\bf 2}(3) (1969), ~311--395.

\bibitem{haase}
 M.~Haase, {\it The Functional calculus for
    sectorial operators},  Operator Theory: Adv. and Appl.
{\bf 169}  Birkhauser Verlag.   Basel-Boston-Berlin, 2006


\bibitem{hie} M. Hieber, S. Monniaux, {\it Heat kernels and
    maximal     $L_{p}-L_{q}$   estimates: The nonautonomous case}, J. Fourier Anal. Appl.
   {\bf 328} (2000) 467-481.

\bibitem{hie1} M. Hieber, S. Monniaux, {\it Pseudo-differential
    operators and maximal regularity results for non-autonomous
    parabolic equations}, Proc. of the AMS, {\bf  128} (1999),
    1047-1053.

\bibitem{kun} P.C. Kunstmann, L. Weis, {\it Maximal $L_{p}$
    regularity    for     parabolic equations, Fourier multiplier theorems and
    $H^{\infty}$ functional calculus}, in: M. Iannelli, R.
    Nagel, S. Piazzera (Eds.), Proceedings of the Autumn School
    on Evolution Equations and Semigroups, in: Levico Lectures,
    {\bf  69}, Springer-Verlag. Heidelberg, 2004,  65-320.


\bibitem{ulm} H. Laasri,  O.E. Agadir {\it Stability for
    non-autonomous linear evolution equations with
    $L_{p}$-maximal regularity}, Czechoslovak Math. J., {\bf 63
    (138)} (2013), 887-908.

\bibitem{laa} H. Laasri, {\it Regularity properties for evolution
    family governed by non-autonomous forms.} 2017.
    https://arxiv.org/abs/1706.06340.

\bibitem{lun}  A. Lunardi, {\it Analytic Semigroups and Optimal
    Regularity     in Parabolic Problems}, Progr. Nonlinear Differential
    Equations Appl., vol. 16, Birkhauser, Basel, 1995.

\bibitem{mey}
 M. Meyries and R. Schnaubelt, Interpolation, embeddings and traces for
anisotropic fractional Sobolev spaces with temporal weights. J.
Funct. Anal., {\bf 262} (2012), 1200-1229.


\bibitem{qua} M. Ouhabaz, {\it Maximal regularity for
    non-autonomous evolution equations governed by forms having
    less regularity}, Arch. der Math. {\bf  105} (2015), 79-91.


\bibitem{pru} J. Pr\"{u}ss, R. Schnaubelt, {\it Solvability
    and     maximal     regularity of parabolic evolution equations with
    coefficients continuous in time}, J. Math. Anal. Appl. {\bf
        256}     (2001), 405-430.

\bibitem{pru1}  J. Pruss,  G. Simonett, {\it
    Maximal regularity for evolution equations in weighted
    $L_p$-spaces}, Arch. Math.,  \textbf{82} (2004), 415--431.

\bibitem{pya1} S. G. Pyatkov,  M. V. Uvarova, {\it Some
    Properties of Solutions of the Cauchy
    Problem for Evolution Equations}, Diff. Equat., {\bf 48} (2012),
    379-389.

    \bibitem{rud} W. Rudin, {\it  Functional analysis}. McGrow-Hill company, New York, 1973.


\bibitem{tan} H. Tanabe,   {\it  Functional Analytic Methods
    for Partial Differential Equations}, Marcel Dekker, Inc.
    New -York, 1997.


\bibitem{trib} H. Triebel, {\it Interpolation Theory, Function
    Spaces, Differential Operators}, North-Holland Mathematical
    Library, {\it 18}, North-Holland Publishing, Amsterdam,
    1978.

\bibitem{tri1}
H. Triebel,  {\it  Theory of Function Spaces} II. Birkhauser Verlag, Basel, 1992.

\bibitem{uva} M.V. Uvarova, {\it On some nonlocal boundary
    value
    problems for evolution equations},  Sib. Adv.  Math.  2011, {\bf
    21}, 211-231.


\bibitem{yagi} A. Yagi,   {\it Abstract Parabolic Evolution
    Equations     and their Applications}, Springer-Verlag. Berlin Heidelberg,
    2010.

\end{thebibliography}
\end{document}